\documentclass{article}

\PassOptionsToPackage{numbers, compress}{natbib}

     \usepackage[final]{neurips_2024}

\usepackage[utf8]{inputenc} 
\usepackage[T1]{fontenc}    
\usepackage{hyperref}       
\usepackage{url}            
\usepackage{booktabs}       
\usepackage{amsfonts}       
\usepackage{nicefrac}       
\usepackage{microtype}      
\usepackage{xcolor}         

\usepackage{amsmath}
\usepackage{amsthm}
\usepackage{graphicx}
\usepackage{caption}
\usepackage{subcaption}
\usepackage{float}
\usepackage{mathtools}
\allowdisplaybreaks

\title{Lee and Seung (2000)'s Algorithms for Non-negative Matrix Factorization: A Supplementary Proof Guide}

\author{%
    Sungjae Cho$^{1,2,3}$\\
    $^1$Department of Computer Science and Operations Research, Université de Montréal\\
    $^2$Architectures of Biological Learning Lab, CHU Sainte-Justine\\
    $^3$Mila -- Quebec Artificial Intelligence Institute\\
    Montréal, Québec, Canada \\
    \texttt{sungjae.cho@umontreal.ca} \\
}

\begin{document}

\maketitle

\begin{abstract}
    Lee and Seung (2000) introduced numerical solutions for non-negative matrix factorization (NMF) using iterative multiplicative update algorithms.
    These algorithms have been actively utilized as dimensionality reduction tools for high-dimensional non-negative data and learning algorithms for artificial neural networks.
    Despite a considerable amount of literature on the applications of the NMF algorithms, detailed explanations about their formulation and derivation are lacking.
    This report provides supplementary details to help understand the formulation and derivation of the proofs as used in the original paper.
\end{abstract}

\section{Introduction}

Lee and Seung (2000)~\cite{Lee2000AlgorithmsFN} introduced numerical solutions for \textbf{\textit{non-negative matrix factorization}} (NMF)~\cite{Gillis2020NonnegativeMF} using iterative multiplicative update algorithms.
These algorithms and their variants~\cite{Gillis2020NonnegativeMF,Wang2013NonnegativeMF,Gan2021NonnegativeMF,SaberiMovahed2024NonnegativeMF,Chen2022ASO,DeHandschutter2021ASO} have been actively utilized as dimensionality reduction tools for high-dimensional non-negative data (e.g., faces, document semantics)~\cite{Lee1999LearningTP} and learning algorithms for artificial neural networks~\cite{Spratling2009UnsupervisedLO}.
Despite a considerable amount of literature on the applications of the NMF algorithms~\cite{Guo2024TheRO}, details about the original paper's formulation and derivation are lacking.
Detailed derivation using matrix algebra elegantly shows the whole picture of the proofs~\cite{Burred2014DetailedDO}, but this work does not directly supplement the original paper.
However, this report provides supplementary details about element-wise proofs, which keeps the description of the original paper. 
This leads to a better understanding of the proofs in the original paper and its NMF algorithms.

\section{Non-negative matrix factorization}

\subsection{NMF problem}

\paragraph{Definition}

NMF in the original paper is the problem of numerically finding non-negative matrix factors $W$ and $H$ given a non-negative matrix $V$ such that $V \approx WH$.

\paragraph{Matrices in practice}

A collection of data vectors is the matrix $V \in \mathbb{R}^{n \times m}$ where $n$ is the number of the dimensions of a single data sample, and $m$ is the number of examples in the data set.
Then, the $i$-th column vector $V_{:i}$ can be a single data sample such as a flattened image vector.
$V (\approx WH)$ is approximately factorized into the matrices $W \in \mathbb{R}^{n \times r}$ and $H \in \mathbb{R}^{r \times m}$.
$W$ and $H$ are usually considered as a transformation function of either $H$ or $W$ and as a collection of factorized vectors.
$r$ is considered as the compressed dimensions of factorized vectors: accordingly, $r \le n$.
In practice, the $i$-th column vector $H_{:i}$ is a compressed vector, and the $j$-th row vector $W_{j:}$ is the coefficients of a linear combination of the $H_{:i}$ elements to approximate the $j$-th element $V_{ji}$ of the $i$-th data column vector $V_{:i}$.
$$
    V_{ji}
    \approx
    \sum_{k} H_{ki} W_{jk}
    =
    H_{:i} W_{j:}
$$

\subsection{Zero elements of $W$ and $H$}

Every element of $W$ and $H$ is non-negative;
however, if an element ($w_a$ or $h_a$) is zero, then it cannot be updated by the suggested multiplicative rules.
Hence, zero elements are redundant, resulting in initializing all elements as positive in practice.
The following proofs are trivial for the zero cases: $w_a = 0$, $h_a = 0$.
Therefore, it is permitted to assume that every element of $W$ and $H$ are positive.

\section{Cost functions}

$A$ and $B$ are non-negative matrices.

\subsection{Euclidean distance}

The cost function to minimize the Euclidean distance is
$$
    C_1(A,B)
    =
    \frac{1}{2} \| A - B \|^2
    =
    \frac{1}{2} \sum_{ij} ( A_{ij} - B_{ij} )^2
    .
$$
This cost function is multiplied by $\frac{1}{2}$ to simplify its derivative.

\subsection{Kullback-Leibler divergence}

Kullback-Leibler (KL) divergence widely defined is
$$
    D_{KL}(A \parallel B) = \sum_{ij} \left ( A_{ij} \log \frac{A_{ij}}{B_{ij}} \right )
    .
$$
However, the KL divergence used as a cost function is
$$
    C_2(A,B)
    =
    D(A \parallel B)
    =
    \sum_{ij} \left ( A_{ij} \log \frac{A_{ij}}{B_{ij}} - A_{ij} + B_{ij} \right )
$$
since $\sum_{ij}A_{ij} = \sum_{ij}B_{ij} = 1$.
The elements of either $A$ or $B$ can be interpreted as being derived from a probability distribution, as they are non-negative and sum to 1.
In practice, $\sum_{ij}A_{ij}$ and $\sum_{ij}B_{ij}$ are not necessarily 1.
$A_{ij}$ and $B_{ij}$ must be strictly positive for the validity of $\log \frac{A_{ij}}{B_{ij}}$ and are recommended to lie within the interval $(0,1]$.
Set $\log 0 = 0$ to avoid implementation errors if the cost function has $\log 0$ terms.

This KL divergence $C_2$ is called the \textbf{\textit{generalized KL divergence}}, which is employed in contexts where the inputs loosely originate from probability distributions.
Hereafter, unless otherwise specified, KL divergence refers to the generalized KL divergence, except in Section~\ref{section:GKLD}.

\subsection{Two alternative formulations of NMF as optimization problems}

\paragraph{Problem 1}
\textit{
    Minimize $\| V - WH \|^2$ with respect to $W$ and $H$, subject to the constraints $W,H \ge 0$.
}

\paragraph{Problem 2}
\textit{
    Minimize $D(V \parallel WH)$ with respect to $W$ and $H$, subject to the constraints $W,H \ge 0$.
}
\\

Although the two cost functions are convex in $W$ only or $H$ only, they are not convex in both variables together (see Appendix~\ref{appendix_nonconvexity} for details about this non-convexity).
Thus, the multiplicative rules suggested in Section~\ref{sec:MU-rules} find local minima rather than global ones.

\section{Multiplicative update rules}
\label{sec:MU-rules}

\paragraph{Theorem 1}

\textit{
    The Euclidean distance $\| V - WH \|$ is nonincreasing under the update rules
}
$$
    H_{a\mu} \leftarrow H_{a\mu} \frac{(W^TV)_{a\mu}}{(W^TWH)_{a\mu}},
    \ \ \ \ \ \
    W_{ia} \leftarrow W_{ia} \frac{(VH^T)_{ia}}{(WHH^T)_{ia}}.
$$
\textit{
    The Euclidean distance is invariant under these updates if and only if $W$ and $H$ are at a stationary point of the distance.
}

\paragraph{Theorem 2}

\textit{
    The divergence $D(V \parallel WH)$ is nonincreasing under the update rules
}
$$
    H_{a\mu} \leftarrow H_{a\mu} \frac{\sum_{i} W_{ia}V_{i\mu} / (WH)_{i\mu}}{\sum_{k} W_{ka}},
    \ \ \ \ \ \
    W_{ia} \leftarrow W_{ia} \frac{\sum_{\mu} H_{a\mu}V_{i\mu} / (WH)_{i\mu}}{\sum_{\nu} H_{a\nu}}.
$$
\textit{
    The divergence is invariant under these updates if and only if $W$ and $H$ are at a stationary point of the divergence.
}

We prove the update rules only for $H$ since those for $W$ can be derived from the following relationships induced by replacing $(V,W,H)$ with $(V^T,H^T,W^T)$.
\begin{align*}
&\min_{H}{\| V - WH \|}
    \xRightarrow{\text{replace}}
    \min_{W^T}{\| V^T - H^T W^T \|}
    =
    \min_{W^T}{\| V^T - (WH)^T \|}
    \Longleftrightarrow
    \min_{W}{\| V - WH \|}
    \\
&\min_{H}{D(V \parallel WH)} 
    \xRightarrow{\text{replace}}
    \min_{W^T}{D(V^T \parallel H^T W^T)}
    =
    \min_{W^T}{D(V^T \parallel (WH)^T)}
    \Longleftrightarrow
    \min_{W}{D(V \parallel WH)}
\end{align*}

The ending parts of Sections~\ref{subsection:gd_squared} and~\ref{subsection:gd_kld} demonstrate that proofs for $H$ can be reused to validate those for $W$.

\section{Multiplicative versus additive update rules}
\label{section:mu}

The multiplicative updates in Theorems 1 and 2 can be derived from the gradient descent formulation.

\subsection{Gradient descent to reduce Euclidean distance}
\label{subsection:gd_squared}

Consider the gradient descent of an element variable $H_{a\mu}$ of $H$.
$$
    H_{a\mu} \leftarrow H_{a\mu} - \eta_{a\mu} \frac{\partial C_1 ( V,  WH )}{\partial H_{a\mu}}
$$
Derive the gradient $\frac{\partial C_1 ( V , WH )}{\partial H_{a\mu}}$.
\begin{align*}
    \frac{\partial C_1 ( V , WH )}{\partial H_{a\mu}}
    &=
    \frac{\partial}{\partial H_{a\mu}} \frac{1}{2} \sum_{ij} ( V_{ij} - (WH)_{ij} )^2 \\
    &=
    \frac{\partial}{\partial H_{a\mu}} \frac{1}{2} \sum_{ij} \left ( V_{ij} - \sum_{k}W_{ik}H_{kj} \right )^2 \\
    &= 
    \frac{1}{2} \cdot 2 
    \sum_{i} \left ( V_{i\mu} - \sum_{k}W_{ik}H_{k\mu} \right )
    \frac{\partial}{\partial H_{a\mu}} \left ( V_{i\mu} - \sum_{k}W_{ik}H_{k\mu} \right ) \\
    &= 
    \frac{1}{2} \cdot 2
    \sum_{i} \left ( V_{i\mu} - \sum_{k}W_{ik}H_{k\mu} \right )
    (-W_{ia})
    \\
    &= 
    - 
    \sum_{i} \left ( W^T_{a i}V_{i\mu} 
    - \sum_{k}W^T_{ai}W_{ik}H_{k\mu} \right )
    \\
    &= 
    - 
    \left [
        \sum_{i} W^T_{ai}V_{i\mu} 
        - \sum_{ik}W^T_{ai}W_{ik}H_{k\mu}
    \right ]
    \\
    &=
    - 
    \left [
        (W^TV)_{a\mu} 
        - (W^TWH)_{a\mu}
    \right ]
\end{align*}
Incorporate this gradient into the gradient descent update rule above.
\begin{align*}
    H_{a\mu}
    \leftarrow& \
    H_{a\mu} + \eta_{a\mu} \left [ (W^TV)_{a\mu} - (W^TWH)_{a\mu} \right ]
    \\
    =& \
    H_{a\mu} - \eta_{a\mu} (W^TWH)_{a\mu} + \eta_{a\mu} (W^TV)_{a\mu}    
\end{align*}
To ensure that the update is multiplicative with respect to $H_{a\mu}$, one possible solution is as follows:
\begin{enumerate}
    \item $\eta_{a\mu}$ must include $H_{a\mu}$ as a factor, allowing all three terms to be factorized by $H_{a\mu}$.
    \item The only negative term $- \eta_{a\mu} (W^TWH)_{a\mu}$, which could make the updated $H_{a\mu}$ negative, is canceled out by the first term $H_{a\mu}$.
\end{enumerate}

The adaptive learning rate 
$$
    \eta_{a\mu}
    =
    \frac{H_{a\mu}}{(W^TWH)_{a\mu}}
$$
satisfies the two conditions above and results in the multiplicative rule as follows.
$$
    H_{a\mu}
    \leftarrow
    \left (
    H_{a\mu}
    - \frac{H_{a\mu}}{(W^TWH)_{a\mu}}(W^TWH)_{a\mu}
    \right )
    + \frac{H_{a\mu}}{(W^TWH)_{a\mu}} (W^TV)_{a\mu}
    =
    H_{a\mu}
    \frac{(W^TV)_{a\mu}}{(W^TWH)_{a\mu}}
$$

The multiplicative rule for $H$ can be derived by replacing $(V,W,H)$ with $(V^T,H^T,W^T)$.
\begin{align*}
    W^T_{a\mu}
    \leftarrow 
    W^T_{a\mu}
    \frac{(HV^T)_{a\mu}}{(HH^TW^T)_{a\mu}}
    &\Longleftrightarrow
    W_{\mu a}
    \leftarrow 
    W_{\mu a}
    \frac{(VH^T)_{\mu a}}{(W H H^T)_{\mu a}}
    \\
    &\Longleftrightarrow
    W_{i a}
    \leftarrow 
    W_{i a}
    \frac{(VH^T)_{i a}}{(W H H^T)_{i a}}    
\end{align*}

\subsection{Gradient descent to reduce KL divergence}
\label{subsection:gd_kld}

Consider the gradient descent of an element variable $H_{a\mu}$ of $H$.
$$
    H_{a\mu} \leftarrow H_{a\mu} - \eta_{a\mu} \frac{\partial C_2 ( V , WH )}{\partial H_{a\mu}}
$$
Derive the gradient $\frac{\partial C_2 ( V , WH )}{\partial H_{a\mu}}$.
\begin{align*}
    \frac{\partial C_2 ( V , WH )}{\partial H_{a\mu}}
    &=
        \frac{\partial}{\partial H_{a\mu}}
        \sum_{ij} \left [ - V_{ij} \log(WH)_{ij} + (WH)_{ij} \right ]
        \\
    &=
        \frac{\partial}{\partial H_{a\mu}}
        \sum_{ij} \left [
            - V_{ij} \log\sum_{k}{(W_{ik}H_{kj})} + \sum_{k}{(W_{ik}H_{kj})}
        \right ]
        \\
    &=
        \frac{\partial}{\partial H_{a\mu}}
        \sum_{i} \left [
            - V_{i\mu} \log\sum_{k}{(W_{ik}H_{k\mu})} + \sum_{k}{(W_{ik}H_{k\mu})}
        \right ]
        \\
    &=
        -
        \sum_{i} \left [
            \frac{V_{i\mu} W_{ia}}{\sum_{k}{(W_{ik}H_{k\mu})}} - W_{ia}
        \right ]
\end{align*}
Put this gradient to the gradient descent update rule above.
\begin{align*}
    H_{a\mu} 
    \leftarrow& \
        H_{a\mu} + \eta_{a\mu}
            \left [
                \sum_{i} \frac{V_{i\mu} W_{ia}}{\sum_{k}{(W_{ik}H_{k\mu})}} - \sum_{i} W_{ia}
            \right ]
    \\
    =& \
        H_{a\mu}
        - \eta_{a\mu} \sum_{i} W_{ia}
        + \eta_{a\mu} \sum_{i} \left [ \frac{V_{i\mu} W_{ia}}{\sum_{k}{(W_{ik}H_{k\mu})}} \right ]
\end{align*}
To ensure that the update is multiplicative with respect to $H_{a\mu}$,
one possible solution is as follows:
\begin{enumerate}
    \item $\eta_{a\mu}$ must include $H_{a\mu}$ as a factor, allowing all three terms to be factorized by $H_{a\mu}$.
    \item The only negative term $- \eta_{a\mu} \sum_{i} W_{ia}$, which could make the updated $H_{a\mu}$ negative, is canceled out by the first term $H_{a\mu}$.
\end{enumerate}

The adaptive learning rate 
$$
    \eta_{a\mu} = \frac{H_{a\mu}}{\sum_{i}W_{ia}}
$$
satisfies the two conditions above and results in the multiplicative rule as follows.
Note that $\sum_{k}{(W_{ik}H_{k\mu})} = (WH)_{i\mu}$.
$$
    H_{a\mu}
    \leftarrow
        \left (
        H_{a\mu}
        - \frac{H_{a\mu}}{\sum_{i}W_{ia}}\sum_{i}W_{ia}
        \right )
        + 
        \frac{H_{a\mu}}{\sum_{i}W_{ia}}
        \sum_{i} \left [ \frac{V_{i\mu} W_{ia}}{\sum_{k}{(W_{ik}H_{k\mu})}} \right ]
    =
        H_{a\mu}
        \frac{\sum_{i} W_{ia}V_{i\mu} / (WH)_{i\mu}}{\sum_{k} W_{ka}}
$$
The multiplicative rule for $H$ can be derived by replacing $(V,W,H)$ with $(V^T,H^T,W^T)$.
\begin{align*}
    W^T_{a\mu}
        \leftarrow 
        W^T_{a\mu}
        \frac{\sum_{i} H^T_{ia}V^T_{i\mu} / (H^TW^T)_{i\mu}}{\sum_{k} H^T_{ka}}
    &\Longleftrightarrow
    W_{\mu a}
        \leftarrow 
        W_{\mu a}
        \frac{\sum_{i} H_{ai}V_{\mu i} / (WH)_{\mu i}}{\sum_{k} H_{ak}}
    \\
    &\Longleftrightarrow
    W_{ia} 
        \leftarrow 
        W_{ia} \frac{\sum_{\mu} H_{a\mu}V_{i\mu} / (WH)_{i\mu}}{\sum_{\nu} H_{a\nu}}    
\end{align*}

\section{Proofs of convergence}

\paragraph{Definition 1}

\textit{
    $G(h,h')$ is an 
    \textbf{auxiliary function}
    for $F(h)$
    if the conditions
    $$G(h,h') \ge F(h), \ \ \ \ \ G(h,h) = F(h)$$
    are satisfied.
}

An objective function that is \textit{difficult to be directly optimized} (i.e., intractable) can be optimized indirectly by optimizing an auxiliary function that is \textit{easier to be directly optimized} (i.e., tractable).
Auxiliary functions facilitate the optimization process.

\paragraph{Lemma 1}

\textit{
    If $G$ is an auxiliary function, then $F$ is nonincreasing under the update
}
$$
    h^{t+1} = \arg \min_{h} G(h,h^t)
    .
$$
\paragraph{\textit{Proof.}}

Prove that $F(h)$ is nonincreasing by showing $F(h^{t+1}) \le F(h^t)$ under the update rule $h^t \rightarrow h^{t+1}$.

\begin{enumerate}
    \item 
        By the definition of auxiliary functions, we have $F(h^{t+1}) \le G(h^{t+1},h^t)$.
    \item 
        By the definition of $h^{t+1}$ in the update rule, $G(h^{t+1},h^t) \le G(h,h^t)$ holds for all $h$, including $h=h^t$.
        Thus, $G(h^{t+1},h^t) \le G(h^t,h^t)$.
    \item Since $G(h,h')$ is an auxiliary function for $F(h)$, we obtain $G(h^t, h^t) \le F(h^t)$.
\end{enumerate}

Therefore, $F(h^{t+1}) \le G(h^{t+1},h^t) \le G(h^t,h^t) \le F(h^t)$, namely, $F(h)$ is nonincreasing under the update rule $h^t \rightarrow h^{t+1}$.

\qed

This iterative algorithm of the Lemma 1 update rule is a \textbf{\textit{majorization-minimization (MM) algorithm}}~\cite{Hunter2004ATO}.
The process of finding the minimum is called \textbf{\textit{minimization}}. 
The conditions for auxiliary functions are called the \textbf{\textit{majorization}} conditions.

\subsection{Convergence of the multiplicative update rules to reduce the Euclidean distance}

\paragraph{Lemma 2}

\textit{
    If $K(h^t)$ is the diagonal matrix (such that its elements are defined as)
}
$$
    K_{ab}(h^t) = \delta_{ab} (W^T W h^t)_{a} / h^t_{a} ,
$$
\textit{
    then
}
$$
    G(h,h^t) = F(h^t) + (h - h^t)^T \nabla F(h^t) + \frac{1}{2} (h - h^t)^T K(h^t) (h - h^t)
$$
\textit{
    is an auxiliary function for
}
$$
    F(h) = \frac{1}{2} \sum_i (v_i - \sum_a W_{ia}h_a)^2 .
$$

\paragraph{\textit{Proof.}}

Note that $W$ is a non-negative matrix,
$h^t$ is a non-negative vector,
and
$\delta_{ab}$ is the Kronecker delta, which filters diagonal elements: $\delta_{ab} = 1$ if $a=b$, and $\delta_{ab} = 0$ otherwise.

To prove that $G(h,h^t)$ is an auxiliary function for $F(h)$, we need to show the two conditions: $G(h,h)=F(h)$ and $G(h,h^t) \ge F(h)$.

\textbf{(i) Prove $G(h,h)=F(h)$:}
If $h^t=h$, then $G(h,h)=F(h) + 0 + 0 = F(h)$.

\textbf{(ii) Prove $G(h,h^t) \ge F(h)$:}
Here, we aim to show $G(h,h^t) - F(h) \ge 0$.
$F(h) 
= \frac{1}{2} \| v - Wh \|^2 
= \frac{1}{2} (v - Wh)^T(v - Wh)$
is a quadratic function of $h$.
Quadratic functions can be rewritten as their second-order Taylor expansion.
Thus, the second-order Taylor expansion at $h=h^t$ of the quadratic function $F(h)$ is exactly $F(h)$.
$$
    F(h) 
    = F(h^t) 
    + (h - h^t)^T \nabla F(h^t) 
    + \frac{1}{2} (h - h^t)^T \nabla^2 F(h^t) (h - h^t)
$$
Compute $G(h,h^t) - F(h)$ by utilizing the second-order derivative $\nabla^2 F(h) = W^TW$ (see Appendix~\ref{appendix_Fh_WTW} for details).
$$
    G(h,h^t) - F(h)
    =
    (h - h^t)^T (K(h^t) - W^T W) (h - h^t)
$$
We show $G(h,h^t) - F(h) \ge 0$ by showing $\nu^T M \nu \ge 0$ for all $\nu (\ne \mathbf{0})$ where the matrix $M$ is defined by $M_{ab} = h^t_a (K(h^t) - W^T W)_{ab} h^t_b$, which is $K(h^t) - W^T W$ element-wisely rescaled by non-negative $h^t_a$ and $h^t_b$.
This holds because $\nu$ covers a broader range than $h - h^t$, and this non-negative rescaling does not change the direction of the inequality.
The primary role of $h^t_a$ and $h^t_b$ is to facilitate achieving $\nu^T M \nu \ge 0$ under the specially defined $K(h^t)$.
\begin{align*}
    \nu^T M \nu \ge 0
    &\Rightarrow
    (h - h^t)^T M (h - h^t) \ge 0
    \\
    &\Leftrightarrow
    (h - h^t)^T h^t_a (K(h^t) - W^T W) h^t_b (h - h^t) \ge 0
    \\
    &\Leftrightarrow
    (h - h^t)^T (K(h^t) - W^T W) (h - h^t) \ge 0
    \\
    &\Leftrightarrow
    G(h,h^t) \ge F(h)
\end{align*}
Let us show $\nu^T M \nu \ge 0$.
\begin{align*}
    \nu^T M \nu
    &= \sum_{ab}   \nu_{a} M_{ab} \nu_{b} \\
    &= \sum_{ab}   \nu_{a} h^t_a (K(h^t) - W^T W)_{ab} h^t_b \nu_{b} \\
    &= \sum_{ab}   \nu_{a} h^t_a K_{ab}(h^t) h^t_b \nu_{b}
       - \sum_{ab} \nu_{a} h^t_a (W^T W)_{ab} h^t_b \nu_{b}
\end{align*}
Simplify the first term $S_1 = \sum_{ab} \nu_{a} h^t_a K_{ab}(h^t) h^t_b \nu_{b}$ as follows.
\begin{align*}
    S_1
    &=
        \sum_{ab} \nu_{a} h^t_a \cdot \delta_{ab}(W^T W h^t)_{a} / h^t_a \cdot  h^t_b \nu_{b} \\
    &=
        \sum_{ab} \nu_{a} \delta_{ab}(W^T W h^t)_{a} h^t_b \nu_{b} \\
    &= 
        \sum_{a} \nu_{a}^2 h^t_a (W^T W h^t)_{a} \\
    &=
        \sum_{a} \nu_{a}^2 h^t_a \sum_{b} (W^T W)_{ab} h^t _{b} \\
    &=
        \sum_{ab} \nu_{a}^2 h^t_a (W^T W)_{ab} h^t _{b}
\end{align*}
Since $W^T W$ is symmetric, namely $(W^T W)_{ab} = (W^T W)_{ba}$, we can rewrite $S_1$ as follows.
\begin{align*}
    S_1
    =
        \frac{1}{2} S_1
        + \frac{1}{2} S_1
    &=
        \frac{1}{2}   \sum_{ab} \nu_{a}^2 h^t_{a} (W^T W)_{ab} h^t _{b}
        + \frac{1}{2} \sum_{ab} \nu_{b}^2 h^t_{b} (W^T W)_{ba} h^t _{a}
        \\
    &=
        \frac{1}{2}   \sum_{ab} \nu_{a}^2 h^t_{a} h^t_{b} (W^T W)_{ab} 
        + \frac{1}{2} \sum_{ab} \nu_{b}^2 h^t_{a} h^t_{b} (W^T W)_{ab}
        \\
    &=
        \sum_{ab} \left ( \frac{1}{2} \nu_{a}^2 + \frac{1}{2} \nu_{b}^2 \right )  h^t_{a} h^t_{b} (W^T W)_{ab}
\end{align*}
Then, we can derive $\nu^T M \nu \ge 0$.
\begin{align*}
    \nu^T M \nu
    &= 
        \sum_{ab} \left ( \frac{1}{2} \nu_{a}^2 + \frac{1}{2} \nu_{b}^2 \right )  h^t_{a} h^t_{b} (W^T W)_{ab}
        - \sum_{ab} \nu_{a} h^t_a (W^T W)_{ab} h^t_b \nu_{b}
        \\
    &= 
        \sum_{ab} \left ( \frac{1}{2} \nu_{a}^2 + \frac{1}{2} \nu_{b}^2 - \nu_a \nu_b \right )  h^t_{a} h^t_{b} (W^T W)_{ab}
        \\
    &= 
        \sum_{ab} \left ( \frac{1}{2} \nu_{a}^2 + \frac{1}{2} \nu_{b}^2 - \nu_a \nu_b \right )  h^t_{a} h^t_{b} (W^T W)_{ab}
        \\
    &= 
        \sum_{ab} \frac{1}{2}  (  \nu_{a} - \nu_{b} )^2  h^t_{a} h^t_{b} (W^T W)_{ab}
        \\
    &\ge
        0
\end{align*}
Hence, $G(h,h^t) \ge F(h)$ is true.

Therefore, $G(h,h^t)$ is an auxiliary function for $F(h)$.

\qed

It appears arbitrary that $(K(h^t) - W^T W)_{ab}$ is multiplied by $h^t_a$ and $h^t_b$ to formulate $M_{ab}$.
However, it is carefully designed for the following reasons:
$h^t_a$ cancels the denominator $h^t_a$ in $K_{ab}$;
$h^t_b$ is replaced by $h^t_a$ due to diagonal symmetry;
finally, this rescaling forms the symmetry of $S_1$ and ensures that $S_1$ is factorized with $h^t_a h^t_b$.

In summary, Lemma 2 confirms that $F(h)$ is nonincreasing under the update rule $h^{t+1} = \arg \min_{h} G(h,h^t)$.
This rule is an instance of the majorization-minimization algorithm that updates $h$ to the minimum of a convex function $G(h,h^{t})$ adjacent to the cost function $F(h)$, which does not increase under this update rule.

\paragraph{Theorem 1}

\textit{
    The Euclidean distance $\| V - WH \|$ is nonincreasing under the update rules
}
$$
    H_{a\mu} \leftarrow H_{a\mu} \frac{(W^TV)_{a\mu}}{(W^TWH)_{a\mu}},
    \ \ \ \ \ \
    W_{ia} \leftarrow W_{ia} \frac{(VH^T)_{ia}}{(WHH^T)_{ia}}.
$$
\textit{
    The Euclidean distance is invariant under these updates if and only if $W$ and $H$ are at a stationary point of the distance.
}

\paragraph{\textit{Proof.}}
Note that $F(h) = \frac{1}{2} \| v - Wh \|^2$ is nonincreasing if and only if $\| v - Wh \|$ is nonincreasing.

The auxiliary function $G$ defined in Lemma 2 is well-designed so that the Lemma 1 update rule (an MM algorithm) yields the multiplicative update rule of Theorem 1.
Let us derive the multiplicative rule of Theorem 1 from the Lemma 1 update rule.

The new update $h^{t+1}$ can be found via the gradient equation $\nabla_h G(h,h^t) = \mathbf{0}$ since $G$ is convex and differentiable.
Derive the $a$-th elements of the gradient and update: $\frac{\partial G(h,h^t)}{\partial h_a}$ and $h^{t+1}_a$.
\begin{align*}
    \frac{\partial G(h,h^t)}{\partial h_a}
    &=
    0 
    + \frac{\partial}{\partial h_a} (h - h^t)^T \nabla F(h^t)
    + \frac{\partial}{\partial h_a} \frac{1}{2} (h - h^t)^T K(h^t) (h - h^t)
    \\
    &=
    \frac{\partial}{\partial h_a} \sum_i (h_i - h^t_i) (\nabla F(h^t))_i
    + \frac{\partial}{\partial h_a} \frac{1}{2} \sum_i (h_i-h^t_i)^2 (W^TWh^t)_i / h^t_i
    \\
    &=
    \frac{\partial}{\partial h_a} (h_a - h^t_a) (\nabla F(h^t))_a
    + \frac{\partial}{\partial h_a} \frac{1}{2} (h_a-h^t_a)^2 (W^TWh^t)_a / h^t_a
    \\
    &=
    (\nabla F(h^t))_a
    + (h_a-h^t_a) (W^TWh^t)_a / h^t_a
\end{align*}
Use $
    (\nabla F(h^t))_a 
    = \frac{\partial C_1}{\partial h_a}
    = - [(W^Tv)_{a} - (W^TWh)_{a}]
$, which is derived from the result of $\frac{\partial C_1}{\partial H_{a\mu}}$.
\begin{align*}
    \frac{\partial G(h,h^t)}{\partial h_a}
    &=
    - (W^Tv)_{a} + (W^TWh^t)_{a}
    + (h_a-h^t_a) (W^TWh^t)_a / h^t_a
    \\
    &=
    - (W^Tv)_{a}
    + \frac{h_a}{h^t_a} (W^TWh^t)_a
    \\
    & = 0
\end{align*}
This equation satisfies $h_a = h^{t+1}_a$.
Solve for $h^{t+1}_a$.
\begin{align*}
    h^t_a (W^Tv)_{a}
    &=
    h_a^{t+1} (W^TWh^t)_a
    \\
    h_a^{t+1}
    &=
    h^t_a\frac{(W^Tv)_{a}}{(W^TWh^t)_a}
\end{align*}
Expand the sample dimension to incorporate $V$ and $H$, and transform this equation as an update rule with an arrow.
$$
    H_{a\mu}^{t+1} = H_{a\mu}^t \frac{(W^TV)_{a\mu}}{(W^TWH^t)_{a\mu}}
    \Longleftrightarrow
    H_{a\mu} \leftarrow H_{a\mu} \frac{(W^TV)_{a\mu}}{(W^TWH)_{a\mu}}
$$
This result demonstrates that the Lemma 1 update rule with the auxiliary function of Lemma 2 matches the update rule of Theorem 1.
Consequently, the update rule of Theorem 1 does not increase the Euclidean distance cost function $C_1=\frac{1}{2}\|V - WH\|^2$ as well as $\|V - WH\|$.
This also holds for the update rule of $W$.
Therefore, Theorem 1 is true.

\qed

\subsection{Convergence of the multiplicative update rules to reduce KL divergence}

\paragraph{Lemma 3}

\textit{
    Define
    \begin{align*}
        G(h, h^t) 
        =& 
        \sum_{i} (v_i \log v_i - v_i) + \sum_{ia} W_{ia}h_a \\
        &
        - \sum_{ia} v_i \frac{W_{ia}h^t_a}{\sum_b W_{ib} h^t_b}
        \left (
            \log W_{ia}h_a - \log \frac{W_{ia}h^t_a}{\sum_b W_{ib} h^t_b}
        \right )
        .
    \end{align*}
    This is an auxiliary function for
    $$
        F(h)
        =
        \sum_i 
        \left [
            v_i \log \frac{v_i}{\sum_a W_{ia}h_a}
            - v_i
            + \sum_a W_{ia} h_a
        \right ]
        .
    $$
    }

\paragraph{\textit{Proof.}}
To prove that $G(h,h^t)$ is an auxiliary function for $F(h)$, we need to show the two conditions: $G(h,h)=F(h)$ and $G(h,h^t) \ge F(h)$.

Rewrite $F(h)$ so that it resembles the terms in $G(h,h^t)$.
\begin{align*}
    F(h)
    &=
    \sum_{i}
        \left [
            v_i \log v_i
            - v_i \log \sum_a W_{ia}h_a
            - v_i
            + \sum_a W_{ia} h_a
        \right ]
    \\
    &=
    \sum_{i}
        \left (
            v_i \log v_i - v_i
        \right )
    + \sum_{ia} W_{ia} h_a
    - \sum_{i} v_i \log \sum_{a} W_{ia}h_a
\end{align*}
Simplify $G(h,h^t) - F(h)$.
\begin{align*}
    G(h,h^t) - F(h)
    &=
        - \sum_{ia} v_i \frac{W_{ia}h^t_a}{\sum_b W_{ib} h^t_b}
            \left (
                \log W_{ia}h_a - \log \frac{W_{ia}h^t_a}{\sum_b W_{ib} h^t_b}
            \right )
        + \sum_{i} v_i \log \sum_{a} W_{ia}h_a
    \\
    &=
        - \sum_{ia} v_i \frac{W_{ia}h^t_a}{\sum_b W_{ib} h^t_b}
            \left (
                \log W_{ia}h_a - \log W_{ia}h^t_a + \log \sum_b W_{ib} h^t_b
            \right )
        + \sum_{i} v_i \log \sum_{a} W_{ia}h_a
\end{align*}

\textbf{(i) Prove $G(h,h)=F(h)$:}
Set $h^t=h$ in $G(h, h^t) - F(h)$.
\begin{align*}
    G(h,h) - F(h)
    &=
        - \sum_{ia} v_i \frac{W_{ia}h_a}{\sum_b W_{ib} h_b}
            \left (
                \log W_{ia}h_a - \log W_{ia}h_a + \log \sum_b W_{ib} h_b
            \right )
        + \sum_{i} v_i \log \sum_{a} W_{ia}h_a
    \\
    &=
        - \sum_{ia} v_i \frac{W_{ia}h_a}{\sum_b W_{ib} h_b} \log \sum_b W_{ib} h_b
        + \sum_{i} v_i \log \sum_{a} W_{ia}h_a
    \\
    &=
        - \sum_{i} v_i \frac{\sum_{a}W_{ia}h_a}{\sum_b W_{ib} h_b} \log \sum_b W_{ib} h_b
        + \sum_{i} v_i \log \sum_{a} W_{ia}h_a
    \\
    &=
        - \sum_{i} v_i \log \sum_b W_{ib} h_b
        + \sum_{i} v_i \log \sum_{a} W_{ia}h_a
    \\
    &= 0
\end{align*}
Hence, $G(h,h) = F(h)$.

\textbf{(ii) Prove $G(h,h^t) \ge F(h)$:}
$-\log x$ is convex where $x \in \mathbb{R}^+$.
Accordingly, it follows Jensen's inequality.
$$
    -\log{\left(\sum_a \alpha_a x_a\right)}
    \le
    \sum_a \alpha_a \left ( - \log \left ( x_a \right ) \right )
$$
where $\alpha_a \in (0,1)$ and $\sum_a \alpha_a = 1$.
Select the special case $x_a = \frac{W_{ia}h_a}{\alpha_a} \in \mathbb{R}^+$.
$$
    -\log{\left(\sum_a W_{ia}h_a \right) 
    \le
    \sum_a \alpha_a \left ( - \log \left ( \frac{W_{ia}h_a}{\alpha_a} \right ) \right )}
$$
Here, take the specific case $\alpha_a = \frac{W_{ia}h^t_a}{\sum_b W_{ib}h^t_b}$, which still satisfies $\alpha_a \in (0,1)$ and $\sum_a \alpha_a = 1$.
\begin{align*}
        -\log{\left(\sum_a W_{ia}h_a \right)} 
    &\le
        \sum_a \frac{W_{ia}h^t_a}{\sum_b W_{ib}h^t_b} \left ( - \log \left ( \frac{W_{ia}h_a\sum_b W_{ib}h^t_b}{W_{ia}h^t_a} \right ) \right )
\end{align*}
Multiply by $v_i (\ge 0)$, and then add $v_i \log v_i - v_i + \sum_a W_{ia}h_a$ on both sides.
\begin{align*}
        - v_i \log{\left(\sum_a W_{ia}h_a \right)}
    &\le
        v_i \sum_a \frac{W_{ia}h^t_a}{\sum_b W_{ib}h^t_b} \left ( - \log \left ( \frac{W_{ia}h_a\sum_b W_{ib}h^t_b}{W_{ia}h^t_a} \right ) \right )
     \\
        - v_i \log{\left(\sum_a W_{ia}h_a \right)}
        + v_i \log v_i - v_i + \sum_a W_{ia}h_a
    &\le
        v_i \sum_a \frac{W_{ia}h^t_a}{\sum_b W_{ib}h^t_b} \left ( - \log \left ( \frac{W_{ia}h_a\sum_b W_{ib}h^t_b}{W_{ia}h^t_a} \right ) \right )
        \\
        &\ \ \ \ \
        + v_i \log v_i - v_i + \sum_a W_{ia}h_a
\end{align*}
Simplify the left and right sides to resemble $F(h)$ and $G(h,h^t)$, respectively.
\begin{align*}
        v_i \log \frac{v_i}{\sum_a W_{ia}h_a}
        - v_i
        + \sum_a W_{ia}h_a
    &\le
        v_i \log v_i - v_i + \sum_a W_{ia}h_a
        \\
        &\ \ \ \ \
        - v_i \sum_a \frac{W_{ia}h^t_a}{\sum_b W_{ib}h^t_b} \left ( \log W_{ia}h_a - \log \frac{W_{ia}h^t_a}{\sum_b W_{ib}h^t_b} \right )
\end{align*}
Sum over $i$ on both sides.
\begin{align*}
        \sum_i 
        \left [
            v_i \log \frac{v_i}{\sum_a W_{ia}h_a}
            - v_i
            + \sum_a W_{ia} h_a
        \right ]
    \le
        &
        \sum_i ( v_i \log v_i - v_i ) + \sum_{ia} W_{ia}h_a
        \\
        &
        - \sum_{ia} v_i \frac{W_{ia}h^t_a}{\sum_b W_{ib}h^t_b} \left ( \log W_{ia}h_a - \log \frac{W_{ia}h^t_a}{\sum_b W_{ib}h^t_b} \right )
\end{align*}
Then, the left and right sides are $F(h)$ and $G(h,h^t)$, respectively.
Thus, $G(h,h^t) \ge F(h)$.

Therefore, $G(h,h^t)$ is an auxiliary function for $F(h)$.

\qed

\paragraph{Theorem 2}

\textit{
    The divergence $D(V \parallel WH)$ is nonincreasing under the update rules
}
$$
    H_{a\mu} \leftarrow H_{a\mu} \frac{\sum_{i} W_{ia}V_{i\mu} / (WH)_{i\mu}}{\sum_{k} W_{ka}},
    \ \ \ \ \ \
    W_{ia} \leftarrow W_{ia} \frac{\sum_{\mu} H_{a\mu}V_{i\mu} / (WH)_{i\mu}}{\sum_{\nu} H_{a\nu}}.
$$
\textit{
    The divergence is invariant under these updates if and only if $W$ and $H$ are at a stationary point of the divergence.
}

\paragraph{\textit{Proof.}}

Like the case of Theorem 1, 
the auxiliary function $G$ defined in Lemma 3 is well-designed so that the Lemma 1 update rule (an MM algorithm) yields the multiplicative update rule of Theorem 2.
Let us derive the multiplicative rule of Theorem 2 from the Lemma 1 update rule.

The new update $h^{t+1}$ can be found via the gradient equation $\nabla_h G(h,h^t) = \mathbf{0}$ since $G$ is convex and differentiable.
Consider the $a$-th elements of the gradient and update: $\frac{\partial G(h,h^t)}{\partial h_a}$ and $h^{t+1}_a$.
\begin{align*}
    \frac{\partial G(h,h^t)}{\partial h_a}
    =
        0 
        + \sum_i W_{ia}
        - \sum_{i} v_i \frac{W_{ia}h^t_a}{\sum_b W_{ib}h^t_b} \frac{1}{h_a}
    =
        0
\end{align*}
This equation satisfies $h_a = h^{t+1}_a$.
Solve for $h^{t+1}_a$.
\begin{align*}
        \sum_i W_{ia}
    &=
        \sum_{i} v_i \frac{W_{ia}h^t_a}{\sum_b W_{ib}h^t_b} \frac{1}{h^{t+1}_a}
    \\
        h^{t+1}_a \sum_i W_{ia}
    &=
        h^t_a \sum_{i} v_i \frac{W_{ia}}{\sum_b W_{ib}h^t_b}
    \\
        h^{t+1}_a
    &=
        h^t_a \frac{\sum_{i} (W_{ia} v_i / \sum_b W_{ib}h^t_b)}{\sum_i W_{ia}}
\end{align*}
Expand the sample dimension to incorporate $V$ and $H$, and transform it as an update rule with an arrow.
\begin{align*}
    &H^{t+1}_{a\mu}
    =
        H^t_{a\mu} \frac{\sum_{i} (W_{ia} V_{i\mu} / \sum_b W_{ib}H^t_{b\mu})}{\sum_i W_{ia}}
    =
        H^t_{a\mu} \frac{\sum_{i} W_{ia} V_{i\mu} / (WH^t)_{i\mu}}{\sum_k W_{ka}}
    \\
    &\Longleftrightarrow
    H_{a\mu}
    \leftarrow
        H_{a\mu} \frac{\sum_{i} W_{ia} V_{i\mu} / (WH)_{i\mu}}{\sum_k W_{ka}}
\end{align*}
This result demonstrates that the Lemma 1 update rule with the Lemma 3 auxiliary function matches the update rule of Theorem 2.
Consequently, the update rule of Theorem 2 does not increase the KL divergence cost function $C_2 = D(V \parallel WH)$.
This also holds for the update rule of $W$.
Therefore, Theorem 2 is true.

\qed

\section{Discussion}

\subsection{Different KL divergence formulation: Generalized KL divergence}
\label{section:GKLD}

$$
    F(h)
    =
    D(v \parallel Wh)
    =
    \sum_i 
    \left [
        v_i \log \frac{v_i}{\sum_a W_{ia}h_a}
        - v_i
        + \sum_a W_{ia} h_a
    \right ]
$$

The original paper does not constrain the elements of $W$ and $H$ to be probabilities, even though they are intended to be used as parameters of KL divergence.
The additional terms $\sum_i ( - v_i + \sum_a W_{ia} h_a )$ to the logarithm terms $\sum_i v_i \log (v_i / \sum_a W_{ia}h_a)$ in $F(h)$ solves this issue.

Here, we employ KL divergence of the scalar $v_i$ rather than of the vector $v$ for simplicity, which can be generalized to the vector $v$ as well as the matrix $V$. 

The original paper does not use KL divergence commonly utilized in statistics, which is defined as 
$$D_{KL}(v_i \parallel W_{i:}h) = v_i \log\frac{v_i}{W_{i:}h}.$$
The NMF algorithms do not force both $v_i$ and $W_{i:}h$ to be probabilities bounded between 0 and 1.
Thus, decreasing KL divergence with respect to $h_a$ does not lead to $v_i \approx W_{i:}h$ but to $h_a \rightarrow \infty$.
Accordingly, KL divergence needs to be reformulated to have the minimum only if $v_i = W_{i:}h$.
The \textbf{\textit{generalized KL divergence}} 
$$D_{GKL}(v_i \parallel W_{i:}h) = v_i \log\frac{v_i}{W_{i:}h} - v_i + W_{i:}h$$ 
has the minimum only if $v_i = W_{i:}h$.

\begin{figure}[H]
    \centering
    \includegraphics[width=0.7\textwidth]{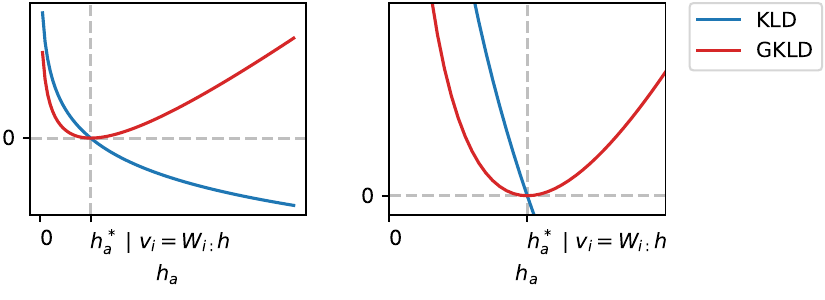}
    \caption{
        The landscapes of KL divergence with respect to $h_a$;
        KLD: KL divergence commonly used;
        GKLD: generalized KL divergence;
        left: larger scale; right: smaller scale.
        }
    \label{fig:GKLD}
\end{figure}

The KL divergence $D_{KL}$ (blue in Fig.~\ref{fig:GKLD}) is a strictly decreasing convex function of $W_{ia}$  since its derivative is negative.
Note that $D_{KL}$ only allows $v_i > 0$ and $W_{ia} > 0$.
\begin{align*}
    \frac{\partial D_{KL}(v_i \parallel W_{i:}h)}{\partial h_a}
    =
        -v_i
        \frac{\partial}{\partial h_a}
        \log \sum_k W_{ik}h_k
    =
        - \frac{v_i W_{ia}}{\sum_k W_{ik}h_k}
    <
        0
\end{align*}
The generalized KL divergence $D_{GKL}$ (red in Fig.~\ref{fig:GKLD}) is a convex function that is neither only increasing nor only decreasing for $W_{ia}$.
This function has the only minimum as the following equation shows.
\begin{align*}
    \frac{\partial D_{GKL}(v_i \parallel W_{i:}h)}{\partial h_a}
    &=
        \frac{\partial}{\partial h_a}
        \left (
            -v_i\log \sum_k W_{ik}h_k
            +
            \sum_k W_{ik}h_k
        \right )
    =
        -\frac{v_i W_{ia}}{\sum_k W_{ik}h_k}
        +
        W_{ia}
    \\
    &=
        \frac{W_{ia}}{\sum_k W_{ik}h_k}
        \left (
            \sum_k W_{ik}h_k - v_i
        \right )
    =
        0
    \\
    &\Leftrightarrow
        v_i
        =
        \sum_k W_{ik}h_k
        =
        W_{i:}h
\end{align*}
$h_a^*$ denotes the element $h_a$, such that the vector $h$ satisfies $v_i = W_{i:}h$.
Keep in mind that $h_a$ is the only dependent variable in this optimization problem.
$D_{GKL}$ demonstrates convexity, exhibiting the following trends across different ranges of $h_a$.
$$
    \begin{cases}
        h_a < h^*_a \ \Rightarrow \ \partial D_{GKL} / \partial h_a < 0 \ \Rightarrow \ \text{$D_{GKL}$ is decreasing.} \\
        h_a = h^*_a \ \Rightarrow \ \partial D_{GKL} / \partial h_a = 0 \ \Rightarrow \ \text{$D_{GKL}$ is at the stationary point.} \\
        h_a > h^*_a \ \Rightarrow \ \partial D_{GKL} / \partial h_a > 0 \ \Rightarrow \ \text{$D_{GKL}$ is increasing.}
    \end{cases}
$$
Hence, this convexity secures the unique minimum for non-probability decision variables and leads to easier optimization by MM algorithms.

In practice, the original paper linearly scales the elements of $V$ so that the pixel mean and standard deviation are equal to $0.25$ and clipped to the interval $[0,1]$~\cite{Lee1999LearningTP}.
Moreover, the sum of the elements of $WH$ is not strictly 1 to be a probability distribution since that of $V$ is not 1 either.

\subsection{Design of the auxiliary functions}

The original paper~\cite{Lee2000AlgorithmsFN} does not elaborate on the specific inspiration behind choosing the special auxiliary functions, likely due to the limited length of the paper.
Nevertheless, it is worth speculating on the motivation for designing these auxiliary functions to understand the update rules better.

The authors appear to have initially induced the suggested multiplicative update operations for $H_{a\mu}$ or $W_{ia}$.
Subsequently, the authors likely designed convex auxiliary functions to produce the update rules.
The authors likely considered multiplicative updates because they are easier to ensure non-negativity than additive/subtractive updates.
Multiplication updates $h \leftarrow h \cdot \Delta h$ preserve non-negativity as long as $h, \Delta h$ are constrained to be non-negative.
In contrast, additive/subtractive updates $h \leftarrow h \pm \Delta h$ ($\Delta h \ge 0$) require careful constraints on the magnitudes of  $|h|$ and $\Delta h$ for decreasing updates $h \leftarrow h - \Delta h$, to ensure non-negativity.

I hypothesize that the auxiliary functions in Lemmas 2 and 3 were designed in the following manner.

\paragraph{Euclidean distance}

The scalar function $F(h)$ is a quadratic function with respect to a vector variable $h$.
A convex quadratic function can be easily derived using the second-order Taylor approximation at $h=h^t$ where $K(h^t)$ implies the curvature at $h=h^t$.
$$
    G(h,h^t) = F(h^t) + (h - h^t)^T \nabla F(h^t) + \frac{1}{2} (h - h^t)^T K(h^t) (h - h^t)
$$
The approximation satisfies $G(h^t,h^t)=F(h^t)$, which is the first condition for being an auxiliary function.

To secure the second condition $G(h^t,h^t) \ge F(h^t)$, the curvature of the approximation $G$ should be not flatter than that of $F$.
Let us find $K$ that can produce the multiplicative update rule under $\frac{\partial G(h,h^t)}{\partial h_a} = 0$.
Assume that $K$ is diagonal for simplicity.
\begin{align*}
    \frac{\partial G(h,h^t)}{\partial h_a}
    &=
        0
        + \frac{\partial}{\partial h_a}
            (h_a - h^t_a) (\nabla F(h^t))_a
        + \frac{\partial}{\partial h_a}
            \frac{1}{2} (h_a - h^t_a)^2 K_{aa}(h^t)
    \\
    &=
        0
        + [ -(W^Tv)_{a} + (W^TWh^t)_{a} ]
        + (h_a-h^t_a) K_{aa}(h^t)
    =
        0
    \\
    \Leftrightarrow
    h_a K_{aa}(h^t)
    &=
    (W^Tv)_{a} - (W^TWh^t)_{a}
    + h_a^t K_{aa}(h^t)
    \\
    \Leftrightarrow
    h_a
    &=
    \frac{(W^Tv)_{a} - (W^TWh^t)_{a}}{K_{aa}(h^t)}
    + h_a^t
    =
    h_a^t
    \left (
        \frac{(W^Tv)_{a} - (W^TWh^t)_{a}}{h_a^tK_{aa}(h^t)}
        + 1
    \right )
\end{align*}
Here, $h_a$ becomes the next $h_a$, namely $h_a^{t+1}$.
The last equation above is equal to the suggested multiplicative update $h_a^{t+1} = h^t_a\frac{(W^Tv)_{a}}{(W^TWh^t)_a}$.
Thus, the following holds.
\begin{align*}
    \frac{(W^Tv)_{a}}{(W^TWh^t)_a}
    &=
    \frac{(W^Tv)_{a} - (W^TWh^t)_{a}}{h_a^tK_{aa}(h^t)}
    + 1
\end{align*}
By subsequently solving for $K_{aa}$, we find that this $K$ matches $K$ defined in Lemma 2.
\begin{align*}
    K_{aa}(h^t)
    &=
    \left ( {(W^Tv)_{a} - (W^TWh^t)_{a}} \right )
    \frac{(W^TWh^t)_a}{h_a^t \left ( (W^Tv)_{a} - (W^TWh^t)_a \right )}
    \\
    &=
    \frac{(W^TWh^t)_a}{h_a^t}
\end{align*}
Now, we know that the curvature of $C_1$ or $F$ is not sharper than that of $K$ since $\frac{\partial^2 C_1}{\partial h_a^2} \le K_{aa}(h^t)$.
The second term $\sum_{i \neq a}(W^TW)_{ai}h^t_i / h_a^t$ of $K_{aa}(h^t)$ below is non-negative due to the non-negativity of $W$ and $h$.
Hence, $\frac{\partial^2 C_1}{\partial h_a^2} \le K_{aa}(h^t)$.
\begin{align*}
    \frac{\partial^2 C_1}{\partial h_a^2}
    &= - \frac{\partial}{\partial h_a}[(W^Tv)_{a} - (W^TWh)_{a}]
    = \frac{\partial}{\partial h_a} \sum_i{(W^TW)_{ai}h_i}
    = (W^TW)_{aa}
    \\
    K_{aa}(h^t)
    &=
    \frac{(W^TWh^t)_a}{h_a^t}
    =
    \frac{\sum_i(W^TW)_{ai}h^t_i}{h_a^t}
    =
    (W^TW)_{aa} + \sum_{i \neq a}(W^TW)_{ai}h^t_i / h_a^t
\end{align*}
As a result, $G$ is located above or equally to $F$: $G(h,h^t) \ge F(h)$.
Thus, $G$ is an auxiliary function for $F$.
Therefore, this special case of $K$ is an auxiliary function that produces the suggested multiplicative update rule.

\paragraph{KL divergence}

The function $F$ in Lemma 3 is 
$$
    F(h)
    =
        \sum_i \left ( v_i \log{v_i} - v_i \right )
        + \sum_{ia} W_{ia} h_a
        - \sum_i v_i \log{\sum_a W_{ia}h_a}
$$
where the only two terms $\sum_{ia} W_{ia} h_a$ and $- \sum_{i} v_i \log \sum_{a} W_{ia}h_a$ are dependent on $h$.
The first $h$-dependent term is linear;
the second one contains negative logarithm terms, which are convex.
Thus, we may transform Jensen's inequality of these logarithms and then obtain a convex auxiliary function.

\textbf{(i) Derive $- \sum_{i} v_i \log \sum_{a} W_{ia}h_a$ from $F(h)$ within Jensen's inequality:}
The negative logarithm $-\log x$ is convex where $x \in \mathbb{R}^+$;
subsequently,
$$
    -\log{\left(\sum_a \alpha_a x_a\right) 
    \le
    \sum_a \alpha_a \left ( - \log \left ( x_a \right ) \right )}
$$
where $\alpha_a \in (0,1)$ and $\sum_a \alpha_a = 1$.
Then, set the special case $x_a = \frac{W_{ia}h_a}{\alpha_a} \in \mathbb{R}^+$ to construct the second $h$-dependent term $- \sum_{i} v_i \log \sum_{a} W_{ia}h_a$ on the left side.
\begin{align*}
    -\log{\left(\sum_a W_{ia}h_a \right)} 
    &\le
    \sum_a \alpha_a \left ( - \log \left ( \frac{W_{ia}h_a}{\alpha_a} \right ) \right )
    \\
    - \sum_i v_i \log{\left(\sum_a W_{ia}h_a \right)} 
    &\le
    - \sum_{ia} v_i \alpha_a \log \left ( \frac{W_{ia}h_a}{\alpha_a} \right )
\end{align*}
\textbf{(ii) Formulate $F(h)$ on the left side:} 
Add the remaining terms of $F$ to both sides.
Then, the left side becomes $F$, and the right side becomes a candidate auxiliary function $g$.
\begin{align*}
    F(h)
    &=
    \sum_i \left ( v_i \log{v_i} - v_i \right )
    + \sum_{ia} W_{ia} h_a
    - \sum_i v_i \log{\left(\sum_a W_{ia}h_a \right)}
    \\
    &\le
    \sum_i \left ( v_i \log{v_i} - v_i \right )
    + \sum_{ia} W_{ia} h_a
    - \sum_{ia} v_i \alpha_a \log \left ( \frac{W_{ia}h_a}{\alpha_a} \right )
    = g(h,h^t)
\end{align*}
\textbf{(iii) Find $\alpha_a$ such that the candidate $g$ satisfies the suggested multiplicative update rule:}
The function $g$ is convex because it is the sum of convex functions, specifically negative logarithms, and the addition of linear functions does not alter its convexity.
Compute the derivative of $g$ with respect to $h_a$.
$$
    \frac{\partial g(h,h^t)}{\partial h_a}
    =
    \sum_i W_{ia}
    - \sum_i v_i \alpha_a \frac{1}{h_a}
$$
Note that the convex function $g$ is not always increasing or decreasing because $\frac{\partial g(h,h^t)}{\partial h_a}$ is not always non-negative or non-positive.
This implies that its extremum is the unique minimum.
Now, the update $h^{t+1} = \arg \min_h g(h,h^t)$ can be found through $\frac{\partial g(h,h^t)}{\partial h_a}=0$.
\begin{align*}
    \frac{\partial g(h,h^t)}{\partial h_a}
    = 0
    \ \
    \Longleftrightarrow
    \ \
    h_a
    =
    h^{t+1}_a
    =
    \frac{\sum_i v_i \alpha_a}{\sum_i W_{ia}}
\end{align*}
The last equation above is equivalent to the suggested multiplicative rule.
\begin{align*}
    h^{t+1}_a 
    = 
    h^t_a \frac{\sum_{i} W_{ia} v_i / \sum_b W_{ib}h^t_b}{\sum_i W_{ia}}
    = 
    h^t_a \frac{\sum_{i} v_i ( W_{ia} / \sum_b W_{ib}h^t_b )}{\sum_i W_{ia}}
\end{align*}
To keep this equivalence, $\alpha_a$ should be 
\begin{align*}
    \alpha_a
    =
    \frac{W_{ia} h^t_a}{\sum_b W_{ib}h^t_b}
    .
\end{align*}
These specially selected $x_a$ and $\alpha_a$ are equal to those given in the proof of Lemma 3.
Hence, $g(h,h^t)$ is exactly the auxiliary function $G(h,h^t)$ introduced in Lemma 3.

\begin{ack}
    This research is supported in part by the FRQNT Strategic Clusters Program (Centre UNIQUE - Centre de rercherche Neuro-IA du Québec) to S.C.
\end{ack}

\bibliographystyle{abbrv}
\bibliography{main}


\newpage

\appendix

\section{Proof of $\nabla^2 F(h) = W^T W$}
\label{appendix_Fh_WTW}
Prove the Hessian $\nabla^2 F(h) = W^T W$ where $F(h) = \frac{1}{2} \sum_i (v_i - \sum_a W_{ia}h_a)^2$.

\paragraph{\textit{Proof.}}

The element of the Hessian $\nabla^2 F(h)$ at $(p,q)$ is $\frac{\partial^2 F}{\partial h_p \partial h_q}$.
\begin{align*}
    \frac{\partial F}{\partial h_q}
    &=
        \sum_i \left (v_i - \sum_a W_{ia}h_a \right )(-W_{iq})
    =
        \sum_i \left ( \sum_a (W_{iq}W_{ia}h_a) - W_{iq}v_i \right )
    \\
    \frac{\partial^2 F}{\partial h_p \partial h_q}
    &=
        \frac{\partial}{\partial h_p}
        \frac{\partial F}{\partial h_q}
    =
        \sum_i W_{iq}W_{ip}
    =
        \sum_i W_{pi}^TW_{iq}
    =
        (W^TW)_{pq}
\end{align*}
Therefore, $\nabla^2 F(h) = W^T W$.

\qed

\section{Proof of the non-convexity with respect to both $W$ and $H$}
\label{appendix_nonconvexity}
To intuitively understand the non-convexity with respect to both $W$ and $H$, let us investigate the convexity with respect to scalar variables $w_a, h_a \in \mathbb{R}$, then row and column vector variables $w, h \in \mathbb{R}^r$, and finally matrix variables $W \in \mathbb{R}^{n \times r}$ and $H \in \mathbb{R}^{r \times m}$.

\subsection{Euclidean distance}

\subsubsection{Prove the convexity for scalars $w_a$ and $h_a$}

Set the cost function $f(w_a,h_a) = (v_a - w_a h_a)^2$.

If $f(\lambda x^{(1)} + (1 - \lambda) x^{(2)}) \le \lambda f(x^{(1)}) + (1 - \lambda) f(x^{(2)})$ for all $x^{(1)} , x^{(2)} \in \{ (w_a , h_a) \}$ and all $\lambda \in [0,1]$, then the function $f$ is convex.
Consider two points $x^{(1)} = (0,0)$ and $x^{(2)} = (\sqrt{v_a},\sqrt{v_a})$;
choose $\lambda = 1/2$.
Then, we obtain the left and right sides ($L$ and $R$) as follows.

$L 
= f(\lambda x^{(1)} + (1 - \lambda) x^{(2)}) 
= f(\frac{1}{2}\sqrt{v_a}, \frac{1}{2}\sqrt{v_a}) 
= (v_a - v_a / 4)^2 
= \frac{9}{16}v_a^2$

$R 
= \lambda f(x^{(1)}) + (1 - \lambda) f(x^{(2)}) 
= \frac{1}{2}f(0,0) + \frac{1}{2}f(\sqrt{v_a}, \sqrt{v_a})
= \frac{1}{2}v_a^2 + 0 
= \frac{1}{2}v_a^2$

$R - L 
= \frac{1}{2}v_a^2 - \frac{9}{16}v_a^2 
= - \frac{1}{16}v_a^2
\le 0$

Therefore, $f(w_a,h_a)$ is not convex.

\qed

\begin{figure}[h!]
    \centering
    \begin{subfigure}[b]{0.24\textwidth}
        \centering
        \includegraphics[width=\textwidth]{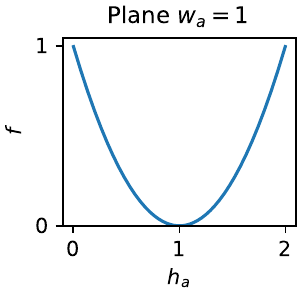}
        \caption{Plane $w_a=1$.}
        \label{fig:w=1}
    \end{subfigure}
    \begin{subfigure}[b]{0.24\textwidth}
        \centering
        \includegraphics[width=\textwidth]{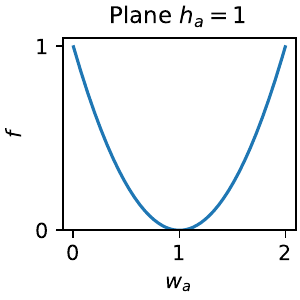}
        \caption{Plane $h_a=1$.}
        \label{fig:h=1}
    \end{subfigure}
    \begin{subfigure}[b]{0.24\textwidth}
        \centering
        \includegraphics[width=\textwidth]{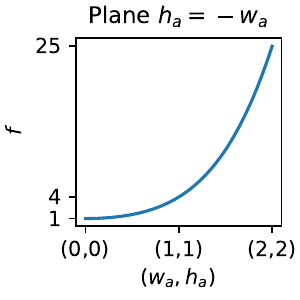}
        \caption{Plane $h_a=-w_a$.}
        \label{fig:h=-w}
    \end{subfigure}
    \begin{subfigure}[b]{0.24\textwidth}
        \centering
        \includegraphics[width=\textwidth]{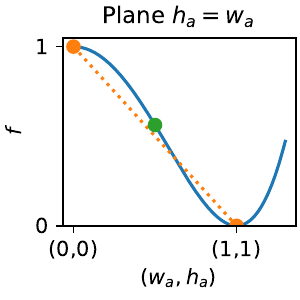}
        \caption{Plane $h_a=w_a$.}
        \label{fig:h=w}
    \end{subfigure}
    \caption{
        The intersection of $f(w_a,h_a) = (v_a-w_a h_a)^2$ where $v_a = 1$ with multiple planes.
        (a,b) These examples show that $f$ is convex if either $w_a$ or $h_a$ is fixed.
        (c) This example shows that $f$ can be convex without fixing $h_a$ or $w_a$.
        (d) This intersection shows a counterexample that proves that $f$ is not always convex. The dots have been used to prove the non-convexity of $f$.
        }
    \label{fig:convexity_SED}
\end{figure}

Figure~\ref{fig:convexity_SED} visualizes what condition $f$ shows convexity or non-convexity.

\subsubsection{Prove the convexity for vectors $w$ and $h$}

Set the cost function $f(w,h) = (v_i - w^Th)^2$
where $w$ is the $i$-th row of $W$, and $h$ is the $i$-th column of $H$.

If $f(\lambda x^{(1)} + (1 - \lambda) x^{(2)}) \le \lambda f(x^{(1)}) + (1 - \lambda) f(x^{(2)})$ for all $x^{(1)} , x^{(2)} \in \{ (w , h) \}$ and all $\lambda \in [0,1]$, then the function $f$ is convex.
Consider two points $x^{(1)} = (\mathbf{0}, \mathbf{0})$ and $x^{(2)} = (s',s')$ where the elements of $s' \in \mathbb{R}^r$ are all $\sqrt{v_i/r}$;
then, choose $\lambda = 1/2$.

$L 
= f(\lambda x^{(1)} + (1 - \lambda) x^{(2)}) 
= f(\frac{1}{2} \cdot \mathbf{0} + \frac{1}{2} s') 
= (v_i - \sum_{i=1}^r \frac{v_i}{4r} )^2 
= (\frac{3v_i}{4} )^2
= \frac{9}{16}v_i^2
$

$R 
= \lambda f(x^{(1)}) + (1 - \lambda) f(x^{(2)}) 
= \frac{1}{2} f(\mathbf{0}, \mathbf{0}) + \frac{1}{2} f(s',s')
= \frac{1}{2}v_i^2 + \frac{1}{2}(v_i - \sum_{i=1}^r \frac{v_i}{r})^2
= \frac{1}{2}v_i^2
$

$
R - L
= \frac{1}{2}v_i^2 - \frac{9}{16}v_i^2
= - \frac{1}{16}v_i^2
\le 0
$

Therefore, $f(w,h)$ is not convex.

\qed

\subsubsection{Prove the convexity for matrices $W$ and $H$}

We have proven the non-convexity of $f$ with respect to a single element $v_a$ of $V$.
The next step is to expand the dimension of a single data sample  from $v_a \in \mathbb{R}, w \in \mathbb{R}^{r}$ to $v \in \mathbb{R}^{n}, W \in \mathbb{R}^{n \times r}$
$$
    v_a \in \mathbb{R} \longrightarrow v \in \mathbb{R}^{n}, \ \ \ \ \
    w \in \mathbb{R}^{r} \longrightarrow W \in \mathbb{R}^{n \times r}
$$
and then expand the dimension of the number of data samples from $v \in \mathbb{R}^{n}, h\in\mathbb{R}^{r}$ to $V\in \mathbb{R}^{n\times m}, H\in\mathbb{R}^{r\times m}$.
$$
    v \in \mathbb{R}^{n}  \longrightarrow V \in \mathbb{R}^{n\times m}, \ \ \ \ \
    h \in \mathbb{R}^{r}  \longrightarrow H \in \mathbb{R}^{r\times m}
$$
These expansions add individual cost terms to the total cost functions $F(W,h)$ and $F(W,H)$ below, with coefficient $1/2$ removed for simplicity.
\begin{align*}
    f(w,h) 
    &= \left [v_i - \sum_{a=1}^{r} w_{a}h_a \right ]^2
    \\
    F(W,h) 
    &= \sum_{i=1}^{n} \left [v_i - \sum_{a=1}^{r} W_{ia}h_a \right ]^2
    \\
    F(W,H) 
    &= \sum_{j=1}^{m} \sum_{i=1}^{n} \left [V_{ij} - \sum_{a=1}^{r} W_{ia}H_{aj} \right ]^2
\end{align*}
A counterexample of convexity can be easily found by reusing the counterexample that proves the non-convexity for vectors $w$ and $h$.
Non-convexity for a pair of vectors $w$ and $h$ can still apply to non-convexity for matrices $W$ and $H$ since the remaining fixed pairs $(w,h)$ are just constant residuals in the total cost function and do not affect the direction of the inequality.
Therefore, the expanded cost function $F$ is not convex with respect to $(W,h)$ as well as $(W,H)$.

\qed

\subsection{KL divergence}

\subsubsection{Prove the convexity for scalars $w_a$ and $h_a$}

Set the cost function $f(w_a,h_a) = v_a \log\frac{v_a}{w_a h_a} - v_a + w_a h_a$.

If $f(\lambda x^{(1)} + (1 - \lambda) x^{(2)}) \le \lambda f(x^{(1)}) + (1 - \lambda) f(x^{(2)})$ for all $x^{(1)} , x^{(2)} \in \{ (w_a , h_a) \}$ and all $\lambda \in [0,1]$, then the function $f$ is convex.
Consider two points $x^{(1)}=(\sqrt{v_a},3\sqrt{v_a})$ and $x^{(2)}=(3\sqrt{v_a},\sqrt{v_a})$;
choose $\lambda = 1/2$.
Then, we obtain the left and right sides ($L$ and $R$) as follows.

$L 
= f(\lambda x^{(1)} + (1 - \lambda) x^{(2)}) 
= f(2\sqrt{v_a}, 2\sqrt{v_a}) 
= v_a \log \frac{v_a}{4 v_a} - v_a + 4 v_a
= v_a (3 - \log{4})$

$R 
= \lambda f(x^{(1)}) + (1 - \lambda) f(x^{(2)}) 
= \frac{1}{2}f(\sqrt{v_a},3\sqrt{v_a}) + \frac{1}{2}f(3\sqrt{v_a},\sqrt{v_a})
\\ 
= 2 \cdot \frac{1}{2} \left ( v_a \log \frac{v_a}{3 v_a} - v_a + 3v_a \right )
= v_a (2 - \log{3})$

$R - L
= v_a (-1 + \log{\frac{4}{3}})
= v_a \log \frac{4}{3e}
\approx v_a \log \frac{4}{3 \cdot 2.718}
= v_a \log \frac{4}{8.154}
< 0$

Therefore, $f(w_a,h_a)$ is not convex.

\qed

\begin{figure}[h!]
    \centering
    \begin{subfigure}[b]{0.24\textwidth}
        \centering
        \includegraphics[width=\textwidth]{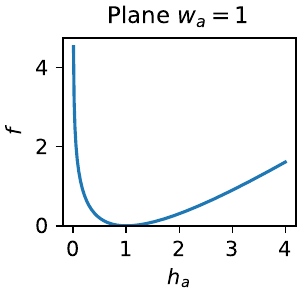}
        \caption{Plane $w_a=1$.}
        \label{fig:w=1}
    \end{subfigure}
    \begin{subfigure}[b]{0.24\textwidth}
        \centering
        \includegraphics[width=\textwidth]{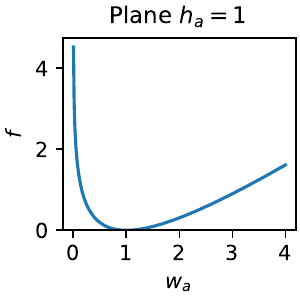}
        \caption{Plane $h_a=1$.}
        \label{fig:h=1}
    \end{subfigure}
    \begin{subfigure}[b]{0.24\textwidth}
        \centering
        \includegraphics[width=\textwidth]{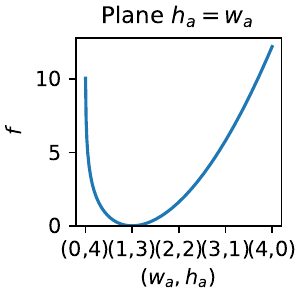}
        \caption{Plane $h_a=w_a$.}
        \label{fig:h=-w}
    \end{subfigure}
    \begin{subfigure}[b]{0.24\textwidth}
        \centering
        \includegraphics[width=\textwidth]{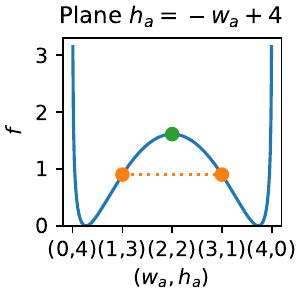}
        \caption{Plane $h_a=-w_a+4$.}
        \label{fig:h=w}
    \end{subfigure}
    \caption{
        The intersection of $f(w_a, h_a) = v_a \log \frac{v_a}{w_a h_a} - v_a + w_a h_a$ where $v_a = 1$ with multiple planes.
        (a,b) These examples show that $f$ is convex if either $w_a$ or $h_a$ is fixed.
        (c) This example shows that $f$ can be convex without fixing $h_a$ or $w_a$.
        (d) This intersection shows a counterexample that proves that $f$ is not always convex. The dots have been used to prove the non-convexity of $f$.
        }
    \label{fig:convexity_KLD}
\end{figure}

Figure~\ref{fig:convexity_KLD} visualizes what condition $f$ shows convexity or non-convexity.

\subsubsection{Prove the convexity for vectors $w$ and $h$}

Set the cost function
$f(w,h) 
    = v_i \log \frac{v_i}{w^T h}
        - v_i
        + w^T h
    = v_i \log \frac{v_i}{\sum_{a=1}^{r} w_a h_a}
        - v_i
        + \sum_{a=1}^{r} w_a h_a
$
where $w$ is the $i$-th row of $W$, and $h$ is the $i$-th column of $H$.

If $f(\lambda x^{(1)} + (1 - \lambda) x^{(2)}) \le \lambda f(x^{(1)}) + (1 - \lambda) f(x^{(2)})$ for all $x^{(1)} , x^{(2)} \in \{ (w , h) \}$ and all $\lambda \in [0,1]$, then the function $f$ is convex.
Consider two points $x^{(1)} = (s', 3s')$ and $x^{(2)} = (3s', s')$ where the elements of $s' \in \mathbb{R}^r$ are all $\sqrt{v_i / r}$;
then, choose $\lambda = 1/2$.
Then, we obtain the left and right sides ($L$ and $R$) as follows.

$L
= f(\lambda x^{(1)} + (1 - \lambda) x^{(2)}) 
= f(2s', 2s')
= v_i \log \frac{v_i}{4 v_i} - v_i + 4 v_i
= v_i (3 - \log{4})$

$R
= \lambda f(x^{(1)}) + (1 - \lambda) f(x^{(2)}) 
= \frac{1}{2}f(s',3s') + \frac{1}{2}f(3s',s')
\\
= 2 \cdot \frac{1}{2} \left ( v_i \log \frac{v_i}{3 v_i} - v_i + 3v_i \right )
= v_i (2 - \log{3})$

$R - L
= v_i (-1 + \log{\frac{4}{3}})
= v_i \log \frac{4}{3e}
\approx v_i \log \frac{4}{3 \cdot 2.718}
= v_i \log \frac{4}{8.154}
< 0$

Therefore, $f(w,h)$ is not convex.

\qed

\subsubsection{Prove the convexity for matrices $W$ and $H$}

$$
    f(w,h) 
    = 
        v_i \log \frac{v_i}{\sum_{a=1}^{r} w_{a} h_a}
        - v_i
        + \sum_{a=1}^{r} w_{a} h_a
$$
$$
    F(W,h) 
    = 
    \sum_{i=1}^{n}
        \left [
            v_i \log \frac{v_i}{\sum_{a=1}^{r} W_{ia} h_a}
            - v_i
            + \sum_{a=1}^{r} W_{ia} h_a
        \right ]
$$
$$
    F(W,H) 
    = 
    \sum_{j=1}^{m}
    \sum_{i=1}^{n}
        \left [
            V_{ij} \log \frac{V_{ij}}{\sum_{a=1}^{r} W_{ia} H_{aj}}
            - V_{ij}
            + \sum_{a=1}^{r} W_{ia} H_{aj}
        \right ]
$$

The same logic of the Euclidean distance case can be applied to this KL divergence case.
Therefore, the cost function $F$ is not convex for $(W,h)$ as well as $(W,H)$.

\qed

\end{document}